\input amstex
\input amsppt.sty
\input epsf
\documentstyle{conm-p}
\NoBlackBoxes

\define\df{\dsize\frac}
\define\bk{\bigskip}

\define \a{\alpha}

\define\tm{\times}

\define\ri{\rightarrow}

\define \s{\sigma}

\define \iy{\infty}

\define\prodl{\prod\limits}

\define \tZ{\tilde{Z}}
\define \tP{\tilde{P}}

\define \la{\langle}
\define \ra{\rangle}

\define \G{\Gamma}

\define \CO{\Cal O}

\define \CP{\Bbb C\Bbb P}
\define \CPN{\Bbb C\Bbb P^N}
\define \CPo{\Bbb C\Bbb  P^1}
\define \CPt{\Bbb C\Bbb P^2}

\define \Z{\Bbb Z}

\define\BP{\Bbb P}

\define \vp{\varphi}

\define \Dl{\Delta}
\define \dl{\delta}
\define \C{\Bbb C}
\define\BC{\Bbb C}

\define \1{^{-1}}
\define \2{^{-2}}
\define \p{\partial}

\define \Aff{\operatorname{Aff}}
\define \Int{\operatorname{Int}}

\define \Gal{\operatorname{Gal}}

\define \Center{\operatorname{Center}}

\define \fc{\frac}

\define \edm{\enddemo}
\define \ep{\endproclaim}

\define\XGal{X_{\Gal}}
\define \CC{{\bbf C}}
\define \PP{{\bbf P}}
\define \bbf{\Bbb}
\topmatter

\title Hirzebruch Surfaces: Degenerations, Related Braid Monodromy, Galois
Covers\endtitle
\author M. Teicher\endauthor
\address Department of Mathematics, Bar-Ilan University, 52900 Ramat-Gan,
Israel\endaddress
\leftheadtext{M. Teicher}
\rightheadtext{Hirzebruch Surfaces: Degenerations, Galois
Covers}

\email teicher\@macs.biu.ac.il\endemail
\thanks   This research was partially supported
by the Emmy Noether Research Institute for  Mathematics, Bar-Ilan
University, and the
Minerva Foundation of Germany.
\endthanks

 \abstract   We describe various properties  of Hirzebruch
surfaces and related constructions: degenerations, braid monodromy,
Galois covers and their Chern numbers.\endabstract
\dedicatory Dedicated to F. Hirzebruch on the occasion of his
70th birthday.\enddedicatory\endtopmatter\document
  \head{\S0.\ Introduction}\endhead

Hirzebruch surfaces were first introduced in 1951, in the paper {\it{``\"Uber
eine
Klasse von einfach-zusammenh\"angenden komplexen Mannigfaltigkeiten''}} (see
\cite{H}). This paper is the first title reprinted in Hirzebruch's
{\it{Gesammelte
Abhandlungen}} (published in 1987 on the occasion of his 60th birthday),
it is
the first part of his dissertation and his very first mathematical paper.
Hirzebruch studied the family of surfaces ~ $\Sigma_n$ for $n \ge 0$ that
are given
by the equation $x_1 y_1^n = x_2   y_2^n$ in $\CC \PP^2 \tm \CC \PP^1.$
He proved that analytically, these 
surfaces are mutually non-isomorphic, whereas
topologically, being $S^2$-bundles over $S^2,$ they fall into only two
homeomorphism
classes, and furthermore, he proved that they are all birationally equivalent.
These surfaces, called {\it{Hirzebruch surfaces}}, have played an important
role
in the theory of algebraic surfaces ever since.
Let us recall the construction as it is usually 
stated nowadays. For $n=k$, the $k$-th  Hirzebruch 
surface is the projectivization of the vector bundle
$\CO_{\CPo}(k)\oplus   \CO_{\CPo}.$ It is usually denoted by $F_k.$
(In fact, any $\CPo$-bundle over $\CPo$ is some $F_k).$

Let $\s$ be a holomorphic section of $\CO_{\CPo}(k)$, and let $E_{0}\subset
F_k$
denote the image of the section $(\s,1)$ of $\CO_{\CPo}(k)\oplus\CO_{\CPo}.$
The curve $E_0$ is called a {\it zero section} of $F_k.$
All zero sections are homologous and hence define a divisor class which is
independent of choice of $\s.$
Let $C$ denote a fiber of $F_k.$
The Picard group of $F_k$ is generated by $E_0$ and $C.$
It is elementary that $E_0^2=k,$ \ $C^2=0$ and $E_0\cdot C=1.$

The surface $F_{0}$ is the quadric $\CC\PP^1 \times \CC\PP^1$, 
and $F_{1}$ is the blow-up of the plane $\CC\PP^2$. For $k>0$, the 
surface~$F_{k}$ contains a unique (irreducible) curve of negative 
self-intersection~$-k$. This curve is a section of the bundle; it is 
denoted~$E_{\infty}$ and it is called the {\it negative section\/} or 
the {\it section at infinity\/}. We mention that it can be contracted 
to an isolated normal singularity, the resulting normal surface being 
the cone over the rational normal curve of degree~$k$. Zero sections 
are always disjoint to~$E_{\infty}$.
Schematically we describe $F_k$ as in Fig. 0.1.


\midinsert
\medskip

\centerline{
\epsfysize=1.5in
\epsfbox{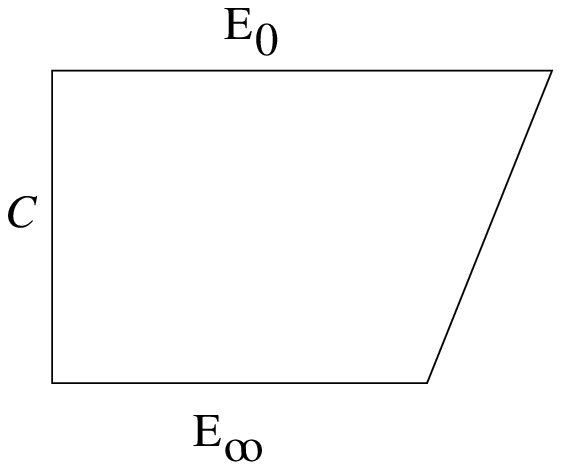}}  

\botcaption{Fig. 0.1}\endcaption
 \endinsert

This might be the place to point out that Hirzebruch, in his mathematical
career,
actually has studied many different classes of surfaces, apart from those
that were
named after him.

In the paper \cite{MoTe1}, published in the year of Hirzebruch's 60th
birthday, we
used the simplest of all Hirzebruch surfaces, namely, the quadric $F_0 =
\BC \BP^1
\tm \BC \BP^1,$ as the starting point to construct a simply connected
surface of
general type with a positive (topological) signature.
That result disproved a famous conjecture in the theory of algebraic
surfaces:  The
{\bf Watershed conjecture} of Bogomolov (see \cite{FH}) stated that a
surface with
non-negative signature should have an infinite fundamental group.
The example was constructed as a Galois cover of $F_0.$
To prove that it is simply connected, its fundamental group was determined by
studying the braid monodromy of the branch curve corresponding to a generic
projection from $F_0,$ suitably embedded in some $\CC\PP^N,$ onto the plane
$\CC\PP^2.$
This work was the starting point of a whole series of papers
\cite{MoTe2} --
\cite{MoTe8}, and \cite{MoRoTe}, \cite{FRoTe}, \cite{Te1} -- \cite{Te4}, in
which we
present our algorithms for computing braid monodromy related to curves,
degeneration
of surfaces, fundamental groups of complements of curves, fundamental groups of
Galois covers of surfaces, and Chern numbers of fibered products.

Some of the examples computed in these papers are based on Galois covers of
Hirzebruch surfaces $F_k.$
In addition to the counterexample, as in \cite{MoTe1}, we produced later
the first
examples of simply connected surfaces of general type with positive
(topological)
signature which are also spin manifolds (\cite{MoRoTe}).
Recall that the signature is positive if $c^2_1/c_2 > 2.$
Corollary 6.3 of this paper gives such an example with $c_1^2/c_2=2.73.$
We also computed an infinite series of pairs of surfaces with the same
Chern numbers
but with different fundamental group, where one group is trivial and the
other of
order going to infinity (\cite{RoTe}).

We believe that fundamental groups of complements of branch curves can
distinguish
among surfaces lying in different connected components of moduli space.
One of our main tools is the braid group (and braid monodromy) technique as
presented in \cite{MoTe3} - \cite{MoTe6}.
The idea to use braid monodromy to compute fundamental groups of complements of
curves started with Van Kampen and Enriques.
Until
the 1980's, very few works dealt with  curves that occur as branch curves
related to
surfaces, in general, and to Hirzebruch surfaces, in particular.
(See sections \S3 and \S4 below for such results).
 One can mention the works \cite{Za} and \cite{Mo}.
It is important to note that the earlier works created a wrong impression
about the
complexity of these fundamental groups, namely, that they 
are ``big", 
and in particular, that they contain free subgroups with two generators.
The braid groups and their close analogues were considered as the typical
examples.
These expectations turned out to be false (see \cite{Te3} for a list of
examples).
The results of Section 3 below are used in \cite{Te5} for the precise
computation of
$\pi_1 (\CP^2 \setminus S),$ where $S$ is the branch curve of a generic
projection of a Hirzebruch surface.

\bigskip

 The paper is divided as follows:

\S0.\ Introduction

\S1.\ Construction and Degeneration of $F_{k(a,b)}$

\S2.\ Braid Monodromy: Definition and Basic Properties

\S3.\ Braid Monodromy Related to Generic Projection of $F_{k(a,b)}$

\S4.\ Fundamental Groups of Galois Covers of Hirzebruch Surfaces

\S5.\ Chern Numbers of Galois Covers

\S6. Intermediate Galois Covers

\bk

\head{\S1.\ Construction and Degeneration of $F_{k(a,b)}$}\endhead

Let $F_k$ be the $k$-th Hirzebruch surface. 
Let $E_0,$\ $E_\iy,$ \ $C$ be as in \S0.
For $a,b\ge1,$ or for $a=0$ and $k\ge1,$
the divisor $aC+bE_0$ on $F_k$ is very ample and thus defines an
embedding $f_{|aC+bE_0|}: F_k\hookrightarrow \CPN.$
Let $F_{k(a,b)}=f_{|aC+bE_0|}(F_k)\ (\subseteq \CPN).$
For $k>0$, the 
map $f_{|0 \cdot C + b E_{0}|}$ collapses the section at infinity to a 
point, so $F_{k(0,b)}$ is the image of the cone over the rational 
normal curve of degree~$k$ with respect to a suitable embedding. 

In \cite{MoRoTe} we constructed a degeneration to a union of $2ab+kb^2$ 
planes in the
following
configuration (in Fig. 1.1, we took $k=2,$\ $a=2,$\ $b=3).$
Each triangle represents a plane and each inner edge represents an
intersection line
between planes.

\midinsert
\centerline{
\epsfysize=1.5in
\epsfbox{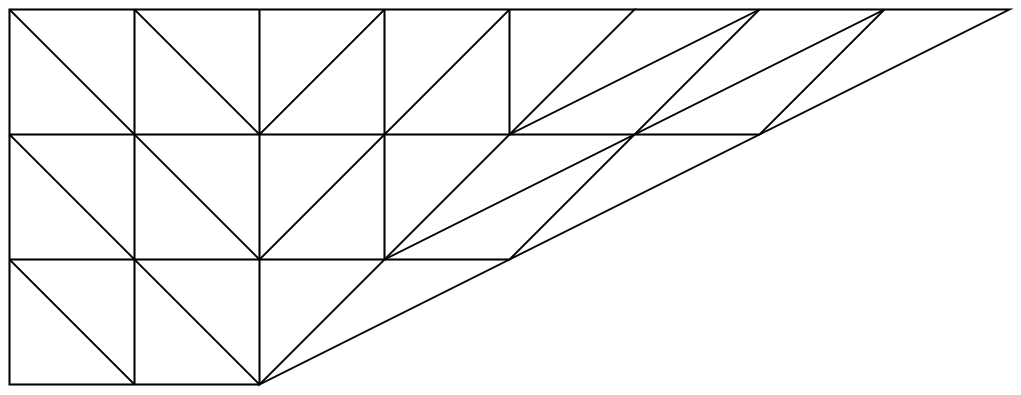}}  

\botcaption{Fig. 1.1}\endcaption
 \endinsert

This degeneration is obtained using a technique developed by us which we
refer to as
the D-construction.
The D-construction is described (and proven to work) in \cite{MoTe5}.
Specific degeneration for the Hirzebruch surfaces using the D-construction is
explained in \cite{MoRoTe}, Section 2 (Theorem 2.1.2).
The difference between the D-construction and other blow-up procedures for
obtaining
degenerations is that we can apply the D-construction also  along a
subvariety of
codim $0$ (see, for example, Step 2 below).
The degeneration is obtained via the following steps:

\roster\item D-construction along $C$ to get $F_{0(1,b)}\cup F_{k(a-1,b)}.$
\item D-construction along $F_{0(1,b)}$ to get $F_{0(1,b)}\cup F_{0(1,b)}\cup
F_{k(a-2,b)}.$
\item Induction  on the second step to get
$\underbrace{F_{0(1,b)}\cup\dots\cup F_{0(1,b)}}_{a\ \text{times}}\cup
F_{k(0,b)}$ (see \cite{MoTe5})\newline

\item Degeneration of each $F_{0(1,b)}$ to a union of $2b$ planes in the
following
configuration (here $b=3)$.

\midinsert
\centerline{
\epsfysize=1.5in
\epsfbox{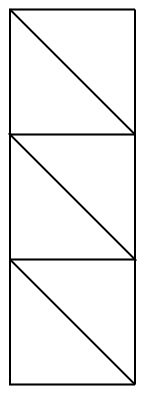}}  

\botcaption{Fig. 1.2}\endcaption
 \endinsert

\item D-construction on $F_{k(0,b)}$ to get
$\underbrace{F_{1(0,b)}\cup\dots\cup
F_{1(0,b)}}_{k\ \text{times}}.$\newline
$(F_{1(0,b)}$ is the Veronese surface $V_b)$
\item Degeneration of each $F_{1(0, b)}$ to a union of $b^2$ planes in the
following
configuration (here $b=3)$:
\midinsert
\centerline{
\epsfysize=1.5in
\epsfbox{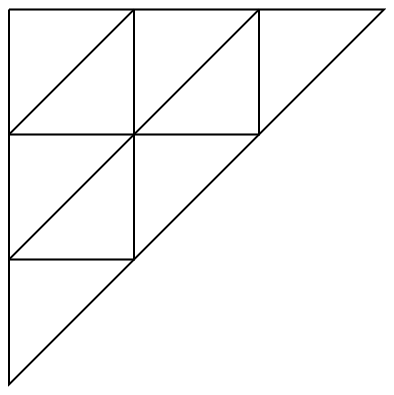}}  
\botcaption{Fig. 1.3}\endcaption
 \endinsert
\endroster

\remark{Remarks}
\roster\item We could go in the ``reverse" direction of the degeneration
and replace
$\underbrace{F_{0(1,b)}\cup\dots\cup F_{0(1,b)}}_{a\ \text{times}} $ by
$F_{0(a,b)}.$
In other words, we might consider a degeneration of $F_{k(a,b)}$ to
  $F_{0(a,b)}\cup \underset k\ \text{times}\to\bigcup V_b.$
\item There are other procedures in progress to obtain a 
degeneration of $V_b$
to a
union of planes.
\item In \cite{CiMiTe} we shall use the above degeneration to describe a new
degeneration of a K3-surface.
\endroster\endremark

\bk
\head{\S 2.\ Braid Monodromy: Definition and Basic Properties}\endhead

In this
section
 we present braid monodromy and braid monodromy
factorizations in
general, and in the next section we shall discuss the one related to Hirzebruch
surfaces.

 Throughout this section (and in section 4)
we shall use the following notations:

$S$ is a curve in $\C^2$ defined over the reals, \ $p=\deg S.$

$\pi: \C^2\ri \C,\  \pi(x,y)=x,$ is the 
first coordinate projection, in a generic coordinate
system defined over the reals.

$K(x)=\{y\bigm| (x,y)\in S\}$.

$N=\{x\bigm| \# K(x)\lvertneqq p\}$ \ (w.l.o.g. $N\subseteq \Bbb R$ since braid
monodromy  is defined up to homotopy type).

$M'=\{(x,y)\in S\bigm| \pi $ is not \'etale at $(x,y)\}\ (\text{clearly},
\pi(M')=N$ and by genericity  $\#(\pi\1(x)\cap M')=1,\ \forall x\in N$).

Let $E$ (resp. $D$)
be a closed disk on the $x$-axis (resp. $y$-axis) such that \linebreak
$M'\subset Int(E\times D),$
($N\subset\Int(E).$)

We choose $u\in\p E,$ real, \quad $x\ll u,\ \ \forall x\in N.$ (Clearly,
$\#(\pi\1(u)\cap S)=p.)$

$K=K(u)=\{q_1,\dots,q_p\}.$

In such a situation,
we are going to introduce {\it braid monodromy.}

\definition{Definition} \ \underbar{Braid monodromy of an affine curve $S$
w.r.t.
$E\times D, \pi, u$}

Every loop  in $E\setminus N$ starting at $u$ has liftings to a system of $p$
paths in \linebreak $(E\setminus N)\times D)\cap S$ starting at $q_1,\dots,
q_p.$
Projecting them horizontally to $D,$ we get $p$ paths
$\{q_1(t),\dots,q_p(t)\}$  in $D$, each one  starts and ends in $K,$
which together can be referred to as a motion.

This motion defines a braid in $B_p[D,K]$ (see \cite{MoTe3}, Section III).
Thus we get a map $\vp:\pi_1 (E\setminus N,u) \ri B_p[D,K].$
This map is evidently a group homomorphism, and it is {\it the braid
monodromy of} $S$ {\it w.r.t.} $E\tm D, \pi, u.$
We sometimes denote\linebreak $\vp$ by $\vp_u.$\enddefinition

It is better to have a notion of braid monodromy of a curve not depending
on the
choice of $D$ and $E,$ when possible and needed:

\definition{Definition} \ $\underline{\text{Braid monodromy of}\  S \
\text{w.r.t.} \ \pi,u}$

Let $\Bbb C_u^1=\{(u,y)\bigm| y\in\Bbb C\}.$ When considering the braid
induced from
the previous motion as an element of the group $B_p [\C_u, K]$ we get the
homomorphism
$\vp:\pi_1(\BC \setminus N, u)\ri B_p[\C_u^1, K]$ which  is called {\it the
braid
mondromy of} $S$ {\it w.r.t.} $\pi, u.$\enddefinition

\midinsert
\centerline{
\epsfysize=2.4in
\epsfbox{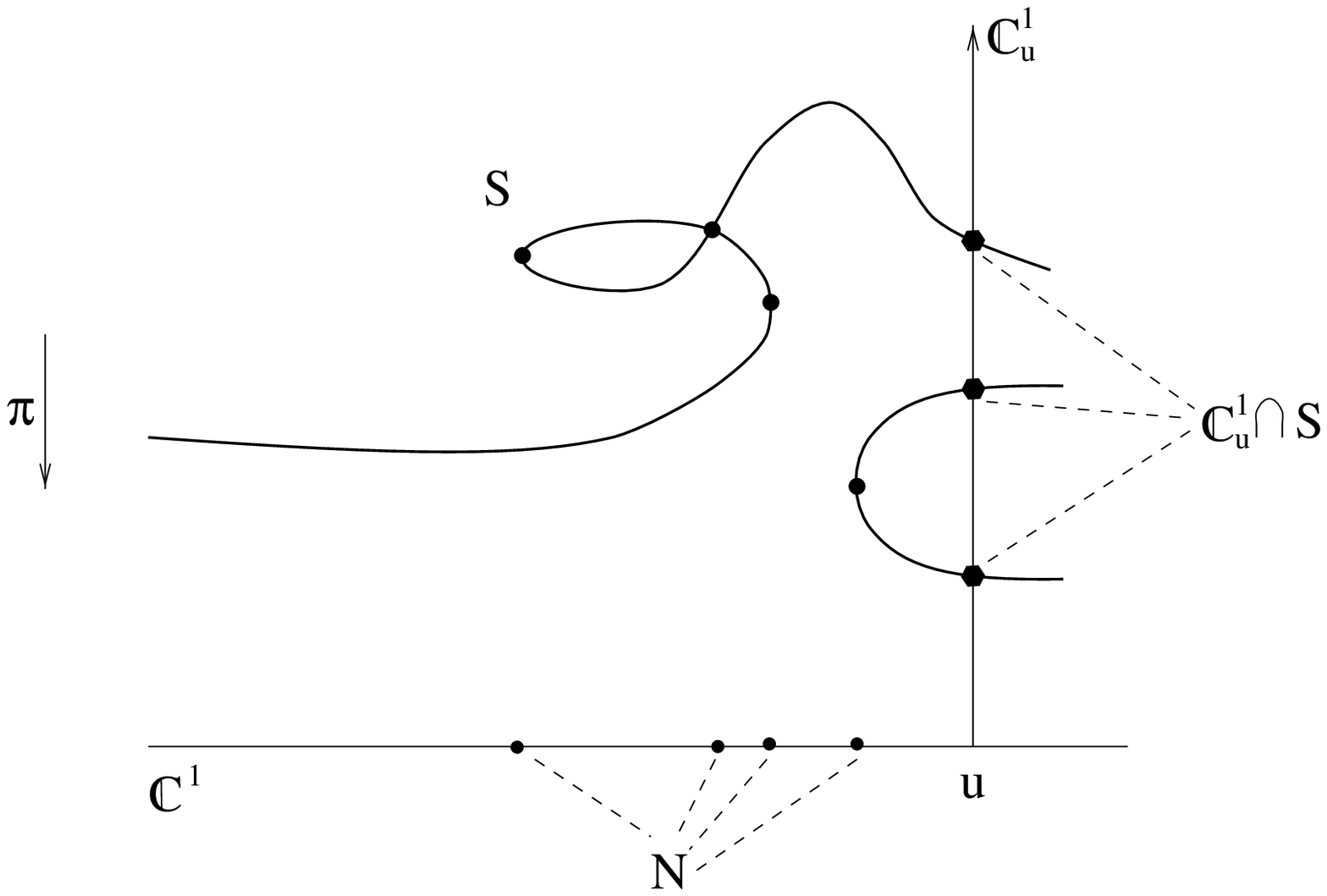}}  
\botcaption{Fig. 2.1}\endcaption
 \endinsert

In order to present an example of a braid monodromy calculation, we recall  a
geometric model of the braid group and the definition of a half-twist.

\definition{Definition}\ $\underline{\text{Braid group}\ B_n[D,K]}$

Let $D$ be
a closed disk in $\Bbb R^2,$ \ $K\subset\Int(D),$ $K$ finite.
Let $B$ be the group of
all diffeomorphisms $\beta$ of $D$ such that $\beta(K) = K\,,\, \beta
|_{\partial D} = \text{Id}\left|_{\partial D}\right.$.
For $\beta_1 ,\beta_2\in B$\,, we
say that $\beta_1$ is equivalent to $\beta_2$ if $\beta_1$ and $\beta_2$ induce
the same automorphism of $\pi_1(D\setminus K,u)$\,.
The quotient of $B$ by this
equivalence relation is called {\it the braid group} $B_n[D,K]$ ($n= \#K$). 
The elements of $B_n[D,K]$ are called {\it braids.} 
We sometimes denote by $\overline\beta$ the braid represented by $\beta.$

\enddefinition

\definition{Definition}\ \underbar{$H(\sigma)$, half-twist defined by
$\sigma$}

 Let
$D,K$ be as above.
Let $a,b\in K,$ and let $\sigma$ be a smooth simple
path in $\Int(D)$ connecting $a$ with $b$ \ s.t. $\sigma\cap
K=\{a,b\}.$
Choose a small regular neighborhood $U$ of $\sigma$ contained in
$\Int(D),$ s.t.
$U\cap K=\{a,b\}$ and an orientation preserving diffeomorphism $f:{\Bbb R}^2
\longrightarrow {\Bbb C}^1$ (${\Bbb C}^1$ is taken with usual ``complex''
orientation) such that
$f(\sigma)=[-1,1]\,,\,$ and $ f(U)=\{z\in{\Bbb C}^1 \,|\,|z|<2\}$\,.
Let $\alpha(r),r\geqslant 0$\,,
be a real smooth monotone function such that $
\alpha(r) = 1$ for $r\in [0,\tsize{3\over 2}]$ and
                $\alpha(r) =   0$ for $ r\geqslant 2.$

Define a diffeomorphism $h:{\Bbb C}^1 \longrightarrow {\Bbb C}^1$ as follows.
For $z\in {\Bbb C}^1\,,\, z= re^{i\varphi},$ let
$h(z) =re^{i(\varphi +\alpha(r)\pi)}$\,.
It is clear that on
$\{z\in{\Bbb C}^1\,|\,|z|\leq\tsize{3\over 2}\}$,\ $h(z)$ is the
positive rotation by $180^{\tsize{\circ}}$ and that
$h(z)=\text{Identity on }\{z\in{\Bbb C}^1\,|\,|z|\ge 2\}$\,, in
particular, on ${\Bbb C}^1 \setminus f(U)$\,.
Considering $(f\circ h\circ f^{-1})|_{D}$ (we always take composition from
left to right), we get a diffeomorphism of $D$ which switches $a$ and $b$
and is
the identity on $D\setminus U$\,.
Thus it defines an element of $B_n[D,K],$ called {\it the
half-twist defined by $\sigma$ and denoted} $H(\sigma).$\enddefinition

The following is the basic braid monodromy associated to a curve with single
singularity.

\proclaim{Proposition - Example  2.1} \ Let $E=\{x\in\C\ |\ |x|\leq
1\},$\ $D=\{y\in\C\ |\ |y| \leq R\},$ $R> 1,$\ $S$ is the curve
$y^2=x^\nu,\ u=1.$
Clearly, here $n=2, N=\{0\},$ \ $ K=\{-1, +1\}$ and $\pi_1(E\setminus N, 1)$ is
generated by
$\G=\p E$ (positive orientation).
Denote by\linebreak $\vp:\pi_1(E\setminus N, 1)\ri B_2[D,K]$ the braid
monodromy of
$S$ w.r.t.
$E\times D, \pi, u.$
Then $\vp(\G)=H^\nu,$ where $H$ is the positive half-twist defined by $[-1,
1]$ (``positive generator" of $B_2[D, K]$).\endproclaim

\demo{Proof} \ We can write $\G = \{e^{2\pi it}, t\in[0,1]\}.$
Lifting $\G$ to $S$ we get two paths:
$$\align
\delta_1(t) &= \left(e^{2\pi it}, \ e^{2\pi i\nu t/2}\right)\\
\delta_2(t) &= \left(e^{2\pi it}, \ -e^{2\pi i\nu t/2}\right).\endalign$$

Projecting $\delta_1(t),$ $\delta_2(t)$ to $D$ we get two paths:
$$\alignat 2
&a_1(t) =e^{\pi it\cdot \nu}, \qquad \qquad  &&0\leq t\leq 1\\
&a_2(t) = -e^{\pi it\cdot \nu}, &&0\leq t\leq 1.
\endalignat$$

This pair of paths $(a_1, a_2),$ each composed of $\nu$ 
consecutive half-circles,  
defines a motion of $ \{1, -1\} $ in $D.$
This motion is the $\nu$-th power of the motion defined by: 
$$\alignat 2
&b_1(t) =e^{\pi it}, \qquad \qquad  &&0\leq t\leq 1\\
&b_2(t) = -e^{\pi it}, &&0\leq t\leq 1.
\endalignat$$

The braid of $B_2\left[D, \{1, -1\}\right]$ induced by this last motion,
 coincides
with
the half-twist $H$ corresponding to $[-1, 1]\subset D.$ Thus $\vp(\G)=H^\nu.$
\hfill$\qed$ \enddemo

\proclaim{Proposition-Definition 2.2} \rom{(Dehn-twist)} Denote by
$ d $ the element of\linebreak
   $\pi_1(D\setminus K,u)$ represented by the loop
$\partial D$ (with positive orientation). There exists a unique element of
$B_n$\,,
denoted by
$\Delta^2_n$ or $\Delta^2_n[D,K]$, such that for any
$\Gamma,$ a simple loop around a single point of $K,$ the (right) action of
$\Dl_n^2$
(as an element of $B_n$) on
$\G$ is as follows: $$ \Gamma\cdot\Delta^2_n = d \Gamma  d ^{-1}\,.
$$
$\Dl_n^2$ is called a Dehn-twist.
\endproclaim
\demo{Proof} \cite{MoTe3}, V.2.1.\edm

\remark{Remark}
Clearly, $\Dl_n^2$ acts as a full-twist around all the points of $K.$
One can justify the notation $\Dl_n^2$, but here we prefer to simply use it
as a
notation.
More about $\Dl_n^2$ can be found in \cite{MoTe3} and \cite{Te6}.\endremark

\proclaim{Proposition 2.3} \rom{(a)}\ $\Dl_n^2\in\Center (B_n).$\qquad
\rom{(b)}\
$\Dl_n^2$ is a product of $n(n-1)$ half-twists.\ep
\demo{Proof} \cite{MoTe3}, V.4.1 and V.2.2.\edm
\proclaim{Proposition - Example 2.4} \ Let $S$ be a union of
$p$ lines, meeting in one point $s_0, s_0=(x_0,y_0).$
Let $D, E, u, K=K(u)$ be as before.
Let $\vp$ be the braid monodromy of $S$ w.r.t. $ E\times D, \pi,u.$
Clearly, here $N=\{x_0\}$ (a single point) and $\pi_1(E\setminus N, u)$ is
generated by
$\G=\p E.$
Then $\vp(\G)=\Dl^2_p=\Dl^2_p \left[ D, K(u)\right].$\endproclaim

\demo{Proof} \ By a continuous change of $s_0$ and of the $n$ lines passing
through
$s_0$ (and by uniqueness of $\Dl^2_p$) we can reduce the proof to the following
case: $S=\bigcup L_k,$  $L_k\:$ $y=j_kx, \ j_k=e^{2\pi ik/p},$ \
$k=0,\dots,p-1.$
Then $N=\{0\}.$
We can take  $u=1,$\ $\G=\{x=e^{2\pi it},$
$t\in[0,1]\},$\quad $K=\{j_k \bigm| k=0\dots p-1\}.$
Lifting $\p E$ to $S$ and then projecting it to $D$, we get $n$ loops:
$$a_k(t)=e^{2\pi i(t+k/p)},\quad k=0,\dots, p-1,\quad t\in[0,1].$$

Thus the motion of the points $a_k(0)=j_k$, represented by the corresponding
loops $a_k(t)$ (for $k=0,\dots,p-1$), is a full-twist which defines
the braid $\Dl^2_p\left[ D,
\{a_k(0)\}\right] =
\Dl^2_p\left[ D, K(1)\right].$
(To check the last fact, see the corresponding actions in
$\pi_1\left(D\setminus K,u\right)$).  \hfill $\qed$\enddemo

\definition{Definition} \ \underbar{Braid monodromy of a projective curve}

Let $B$ be an algebraic curve of degree $p$ in $\C \Bbb P^2.$
Choose generically a line $L$ at infinity  $(\#(L\cap B) = p)$ and affine
coordinates $(x,y)$ in $\C^2 = \C \Bbb P^2 \setminus L$ so that the
 coordinate projection, 
$\pi: \BC^2\to\BC,\ \pi(x,y)=x,$ induces a generic map from 
$B\cap\C^2$ to $\C$
by restriction (in
particular,
the center of this projection in $\C\Bbb P^2$ must lie outside of
$B$). 
Let $N=\{x\in\C \bigm| \pi^{-1}(x)\cap
B\lvertneqq p\},$ $E$ be a closed disk on the
$x$-axis with $N\subset\Int (E) , D$
be a sufficiently large closed disk on the $y$-axis s.t.
$\pi^{-1} (E)\cap B\subset E\times D.$
Choose $u\in\p E.$
Denote by $S = B\cap (E\times D).$
{\it The  braid monodromy of} $B$
{\it w.r.t.} $L, u$ is the braid monodromy of $S$ w.r.t.
$E\times D, \pi, u,$ i.e., the homomorphism
$$\vp:\pi_1(E\setminus N, u)\ri B_p [ D, K]$$\enddefinition
 We recall the notion of a geometric free base of the fundamental group of a
punctured disk in $\BC$ and a basic property of it.
Since we shall choose such
bases both for the $x$-axis and the $y$-axis, we make
independent notations.

\definition{Definition}\ \underbar{A bush}

Let $U$ be a closed disk in $\BC$ and $F$ a finite set in $Int(U),$\
$F=\{w_1,\dots,w_n\},$ \
$v\in\partial U.$

 Consider in $U$ an ordered set of simple
paths
$(T_1,\dots ,T_n)$ connecting the points $w_1,\dots,w_n$ with $v$ such that
\roster
\item $T_i\cap T_j=\{v\}$ if $i\ne j$\,;
\item Each path $T_i$ intersects a small circle around $v$ in a single point
$u_i',$ and the order of these points on the circle is given by the positive
(``counterclockwise'') orientation.
\endroster
We say that two such sets $(T_1,\dots ,T_n)$ and $(T^{\prime}_1,\dots
,T^{\prime}_n)$,  are equivalent if on the homotopy class level we have
$$\ell(T_i)=\ell(T^{\prime}_i)\quad\quad(\text{for }  i=1,\dots ,n)$$
where $\ell(T_i)$ is a closed
loop based at $v$, then following the path $T,$ then encircles $w_i$
counterclockwise and returns (see Fig. 2.2). An equivalence class of such
sets is
called {\it a bush} in
$(U\setminus F,v)$\,. The bush represented by $(T_1,\dots ,T_n)$ is denoted by
$\langle T_1,\dots ,T_n\rangle$.
\enddefinition

\midinsert
\centerline{
\epsfysize=1.5in
\epsfbox{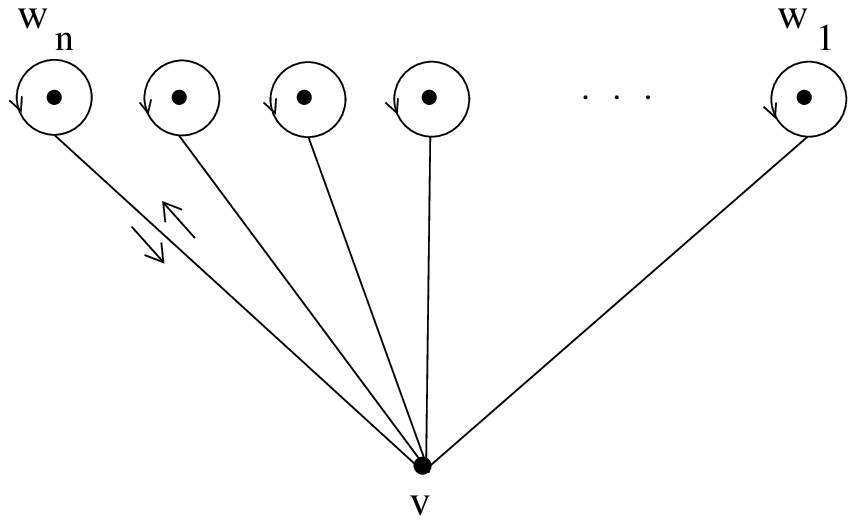}}  
\botcaption{Fig. 2.2}\endcaption
 \endinsert

\newpage

\definition{Definition}\ $\underline{\text{geometric base,}\ g-\text{base}}$

 Let $U,$ $F,$ $v,$   be as above. A $g$-{\it base of} $\pi_1(U\setminus
F,v)$ is an
ordered free base of $\pi_1(U\setminus F,v)$ which has the form
$(\ell(T_1),\dots
,\ell(T_n))$ where $\langle T_1,\dots ,T_n\rangle$ is a bush in $U\setminus
F$ (see
fig. 2.2).
\enddefinition

\proclaim{Proposition 2.5} \ Let $B$ be an algebraic curve of degree $p$ in
$\Bbb C\Bbb P^2.$
Let \newline
$L, \pi,u,D,E,K(u)$ be as in the beginning of \S2.
Let $\vp$ be the braid monodromy of $B$ w.r.t. $L,\pi,u$. Let
$(\delta_1,\dots,\delta_r)$ be a $g$-base
of $\pi_1(E\setminus N,u)$\ $(r=\#N).$ Then
$$\prod^r_{i=1} \vp(\delta_i)=\Dl^2_p=\Dl^2_p [u\times D, K(u) \cap
B].$$\endproclaim

\demo{Proof} \ One can see that $\prod\limits^r_{i=1} \ \delta_i=\p E$\ is
positively oriented (see also \cite{MoTe3}, \S2). Thus we have to prove
that $\vp(\p
E)=\Dl^2_p.$ We can assume $E$ arbitrarily big, so that $\p E$ will be very
close to
$\infty$
at the $x$-axis.
Continuously deforming coefficients in the equations of $B$ such that new
resulting curves always remain transversal to $L,$ we can reduce the proof
to the
case where $B$ is a union of $n$ lines intersecting at a single point. Now use
Proposition - Example 2.4. \hfill $\qed$
\enddemo

Following Proposition 2.5 we define:

\definition{Definition} \
$\underline{\text{Braid monodromy factorization of}\ \Dl_p^2}$ (associated
to a plane
projective
curve)

{\it Braid monodromy factorization of $\Dl_p^2$ (associated to a plane
projective
curve)} is a product  of the form $\Dl^2_p=\prod\limits_i\vp (\delta_i),$ where
$\vp$ is the braid monodromy of the projective curve and $\{\dl_i\}$ is
a $g$-base of $\pi_1(E\setminus N,u).$\enddefinition

\remark{Remarks}

(1)\ A braid monodromy factorization depends, in fact, not only on the curve
but also on the choice of the base. When needed, we then refer to
 {\it braid monodromy factorization of} $\Dl_p^2$
 {\it associated to a curve and a base}\
$\{\dl_i\}.$

(2) \ In the other direction, a $g$-base of $\pi_1(E\setminus N,u)$ and the
corresponding factorization determine the braid
monodromy.
(The values of a homomorphism on a base determine the homomorphism.)
For applications, it is usually sufficient to know a a factorization,
without referencing
to a particular $g$-base (like in the proof of 3.2 below or in  [MoTe7] or 
in [Te5]).

(3) \
For a nonsingular $B,$ each $\vp(\dl_i)$ is a
(positive) half-twist
in $B_p.$ (See \cite{MoTe3}, Prop. IV.1.1).
The associated factorization is then called {\it prime}.

(4) \
A braid monodromy factorization is a
presentation of $\Dl_p^2$ as a product of (positive)
elements in $B_p^+.$
Not all factorizations of  $\Dl_p^2$ to products of positive elements are
induced from a curve. (The definition of positive elements and the semigroup
$B_p^+$ can be found in \cite{MoTe3}, \S5; roughly, these are products of
positive half-twists).\endremark

\proclaim{Proposition 2.6} \ Let $B$ be a
cuspidal curve in $\C \Bbb P^2$ (that is, all
singularities of $B$ are locally of the form
$y^2=x^2$ (a node) or $y^2=x^3$ (a cusp)).
Then any braid monodromy factorization of $\Dl_p^2$ (associated to $B),$
can be written as  a product $\Dl^2_p=\prod\limits_i
(Q_i^{-1} H^{\nu_i}_1 Q_i)$ of suitable conjugates of some fixed
positive half-twist $H_1,$ raised to some power  $\nu_i=1,2, \
\text{or} \ 3.$\endproclaim

\demo{Proof} \ Recall that we are using generic
projections of $\C^2\overset\pi\to\rightarrow \C$ w.r.t. the
projective curve.
Each singularity of $\pi|_B$ is of the type
$y^2=x^\nu,$ $\nu=1,2,$ or $3.$
Now use Proposition - Example 2.1 to get $\vp(\dl_i) = H_i^{\nu_i},$ with $
\nu_i=1,2,$ or
$3,$\quad and where
$H_i$ is a half-twist.
Every two half-twists in $B_p$ are conjugate, so  for every $ i ,$ there
exists a
$Q_i$ s.t.
$H_i=Q_i\1 H_1 Q_i.$
Thus, $\Dl^2_p=\Pi \vp(\dl_i) = \Pi (Q_i\1 H_1^{\nu_i} Q_i).$
\hfill $\qed$
\enddemo

\remark{Remark}  We can take any half-twist for $H_1$.
\endremark

In the next section we shall consider a braid monodromy factorization
related to
Hirzebruch surfaces.

\bk

\head{\S3. \ Braid Monodromy Related to a Generic Projection of
$F_{k(a,b)}$}\endhead

Let $S_{k(a,b)}$ be the branch curve of a generic projection of $F_{k(a,b)}$ to
$\CPt.$
We want to compute the braid monodromy of $S_{k(a,b)}.$
We believe that the ``braid monodromy type" of a branch curve determines the
``deformation type" of the related surface.
Thus our main goal in computing the braid monodromy of a branch curve is to
distinguish between  surfaces which are not a deformation of each other (see
\cite{Te3}). Since $F_{k(a,b)}$ can be deformed to $F_{k-2(a',b')}$ (see
\cite{FRoTe}), it is enough to consider the case
$k=0$ and $k=1.$
The case $k=0$ was presented in \cite{MoTe1}; the case $k=1$ will be
described here.
Theorem 3.2 gives a braid monodromy factorization for $S_{1(a,b)},$ and thus
determine
the braid monodromy type of $S_{1(a,b)}.$
The nonspecialist might want to skip the details of this theorem, and the
subsequent
explanation while realizing that we heavily use the degeneration from \S1
in the
calculation.

Before we state Theorem 3.2, we describe in greater detail the branch curve
of the
degenerated object. We shall use the degeneration of
$F_{1(a,b)}$ described in
\S1.
Recall that  $F_{1(a,b)}$ is degenerated to   $F_{1(a,b)}^0,$ a union of
planes in
the following configuration
(here $b=5$,\ $a=4$):

\midinsert
\centerline{
\epsfysize=2in
\epsfbox{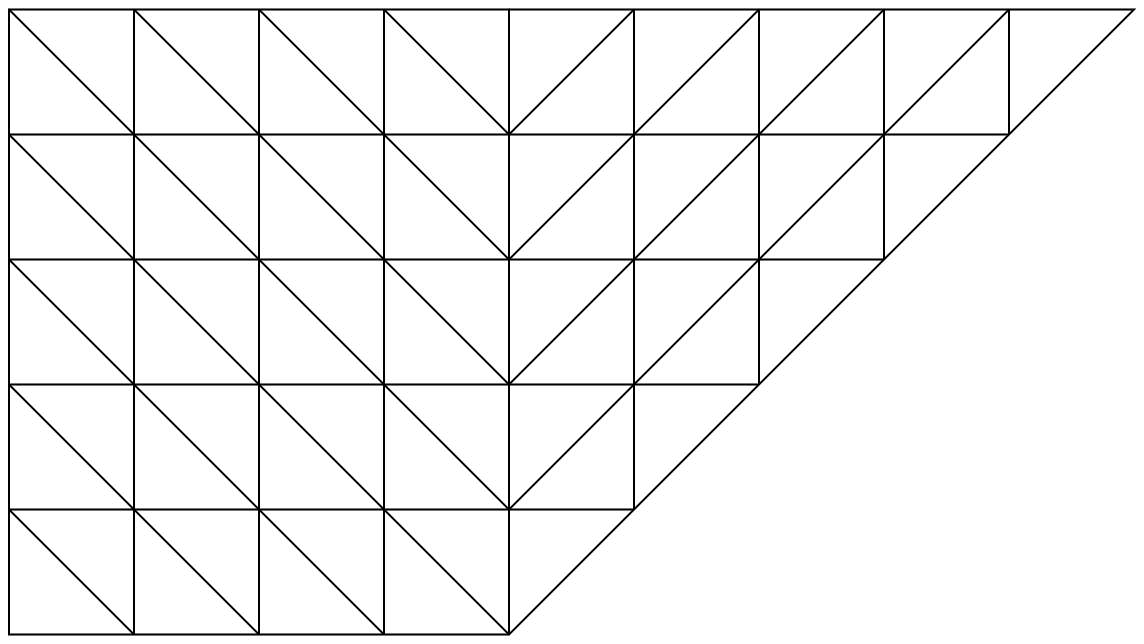}}  
\botcaption{Fig. 3.1}\endcaption
 \endinsert

Each triangle represents a plane and each inner edge represents an
intersection line
between planes.
The number of planes is $2ab+b^2$ and the number of intersection lines is
$3ab-a+\fc{3b}{2}(b-1)$.
We take a generic projection of
$F_{1(a,b)}^0$ onto
$\CPt$ where each plane is projected onto $\CPt.$
The ramification curve of this projection is the union of lines.
The singular points of the ramification curve are represented by 
vertices.
The branch curve of $F_{1(a,b)}^0\to\CPt,$ denoted $S_{1(a,b)}^0,$ is the
image of
the union of lines and its singular points are the images of the vertices
and the
intersection points in $\CP^2$ of the images of any two of the intersection
lines.

 \medskip

We numerate the vertices $a_1,\dots, a _{\nu_{0}}$ from right to left, from
bottom to top
(including two points $a_1$ and $a_{\nu_{0}-b}$ which are not on $S_0$ and  two
points
$a_{m_{0}+b}
$ and
$a_{\nu_{0}}$ which are on $S_0$ but not singular points of $S_0)$, where
$m_0=\fc{b(b+1)}{2}+1$ and  $\nu_0=m_0+a(b+1)+b$
$(=\fc{b(b+1)}{2}+(a+1)(b+1))$.
We numerate the lines in lexicographic order from the bigger index to the
smaller
one: \ $L_1,\dots, L_{p_{0}},$ where
 $p_0=\fc{1}{2}(6ab-2a-3b+3b^2).$ See Fig. 3.2.

\midinsert
\centerline{
\epsfysize=3.4in
\epsfbox{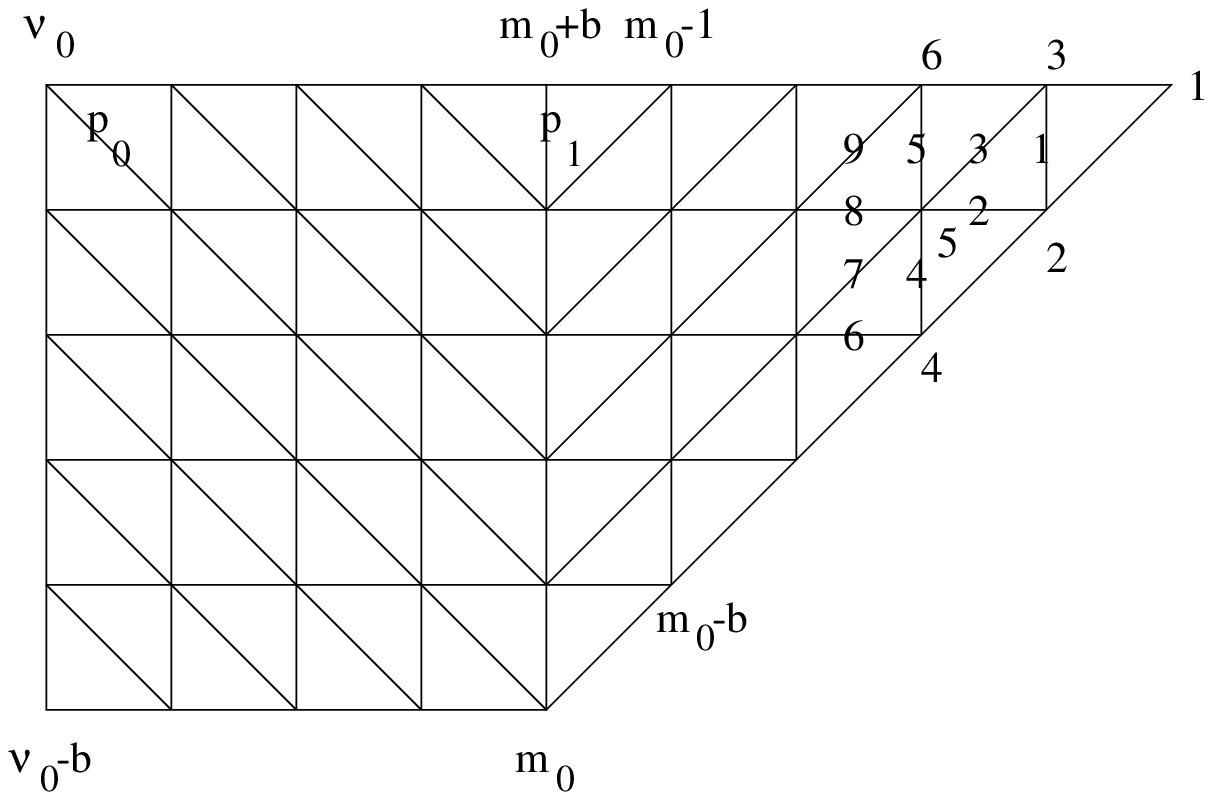}}  
\botcaption{Fig. 3.2}\endcaption
 \endinsert

\proclaim{Lemma 3.1}
$p=\deg S_{1(a,b)}=6ab-2a-3b+3b^2.$\ep
\demo{Proof} Lemma 7.1.3(b) in \cite{MoRoTe} for $k=1$.\edm

\remark{Remark}  In \cite{MoTe4}, \S2,\S3 (see also \cite{MoTe6}, \S1), we
introduced a {\it regeneration process} for ``reconstructing" branch curves
from
the branch
curve of the degenerated object. Since lines are doubled during the
regeneration
process,
$\deg S_{1(a,b)}=2\deg S_{1(a,b)}^0.$
So one can also get the Lemma by doubling the number of intersection lines
in the
above configuration.\endremark

Let $\ell_j$ and $b_i$ be the images of $L_j$ and $a_i$ in $\CPt$ and
$q_i=\ell_i\cap
\C_{u}^1$ (see
\S2).

\proclaim{Theorem 3.2}
The braid monodromy factorization of $\Dl_p^2$ associated to $S_{1(a,b)}$ 
(where $p=\deg S_{1(a,b)}$) is as
follows:\
$\Dl_p^2=\prodl_{i=1}^{\nu_0}\tilde C_i\tilde P_i$ ;
$\tilde P_i$  is the local braid monodromy factorization around $b_i;$ and
$\tilde C_i=  \prod \tZ_{ii',jj'}^2$ for $i<j,$\ $L_i\cap L_j=\emptyset$,
and $\tilde
Z_{i,j}$ is a half-twist
from
$q_i$ to $q_j$ corresponding to the path described in Fig. \rom{3.3},
$j_0$ is the smallest index s.t.    $L_{j_{0}}$
meets
$L_j$ in the vertex 
with the higer index, and $\tilde Z_{ii',jj'}^2=\tilde Z_{ij}^2\tilde
Z_{ij'}^2\tilde Z_{i'j}^2\tilde Z_{i'j'}^2$ (a product of 4 full-twists).
\ep

\midinsert
\centerline{
\epsfysize=0.7in
\epsfbox{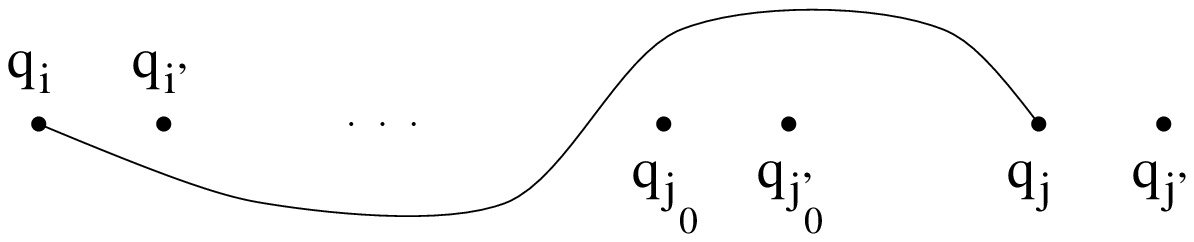}}  
\botcaption{Fig. 3.3}\endcaption
 \endinsert

\demo{Proof} We use the notation of \S2.
Recall that $F_{1(a,b)}$ was degenerated to $F_{1(a,b)}^0$ with the
branch curve $S_{1(a,b)}$ of the generic projection onto $\CP^2$ being
degenerated to $S_{1(a,b)}^0,$ the branch curve of
$F_{1(a,b)}^0\to
\CPt.$ As explained earlier, $S_{1(a,b)}^0$ is an arrangement of  $p_0$ lines.
In fact, $p_0=p/2$ (see the above remark). In this arrangement, 
 no 3 vertices of higher multiplicity (where 3 lines or more meet) 
are collinear.
 In \cite{MoTe3}, \S9,
we computed the   braid monodromy factorization 
associated to such line arrangements.  
It can be
presented as:
$$\Dl_{p_{0}}^2=\prodl_{i=1}^{\nu_0} C_i\Dl_{k_{i}}^2[\C_{u_{0}}^1,S_i\cap
C_{u_{0}}^1]$$
$$C_i= \prod_{i<j, \ L_i\cap L_j=\emptyset} \tZ_{ij}^2$$

\flushpar where $S_i=\{\text{lines through}\ a_i\},$ with
$k_i=\# S_i,$ 
and $\Dl_{k_{i}}^2=\text{Dehn-twist around}$\linebreak $S_i\cap\C_{u_{0}}^1.$
(The notation $\tilde Z_{ij}$ is explained in the formulation of the Theorem.)

By Proposition-Example 2.4, $\Dl_{k_{i}}^2$ is, in fact, the local braid
monodromy
of $S_{1(a,b)}^0$ around $b_i$\ $(\Dl_i^2=1$ for $i=1,m_0+b,\nu_0-b$ and
$\nu_0).$

We apply the regeneration process (from \cite{MoTe4}) on
$\prod\limits_{i=1}^{\nu_0}C_i\Dl_i^2$ to get $\prod\limits_{i=1}^{\nu_0}\tilde
C_i\tilde P_i.$ (See \cite{MoTe4}, \S2 for the starting situation and \S3
for the
regeneration rules I, II, and III (Lemmas 3.1, 3.2 and 3.3)).
When ``regenerating," lines are ``doubled" and each point $q_j$ in the
typical fiber
is replaced by two points $q_j$ and $q_{j'}.$
The product $\tilde C_i$ is easy to describe.
It is the result of applying the second
regenerating rule on
$C_i$, i.e., each full-twist $\tilde
Z_{ij}^2$ is
replaced by the product of the  4 full-twists $\tilde Z_{ij}^2,\tilde
Z_{ij'}^2,\tilde Z_{i'j}^2$,
$\tilde Z_{i'j'}^2$ in this order (written in short
$Z_{ii',jj'}^2)$ (see also C-table in \cite{MoTe1}).

For $a_i$'s which are  not singular points of $S_{1(a,b)}^0,$ we get the
following expression: 
$\tilde
P_{m_{0}+b}=Z_{p_{1}p_{1}'}\ \left(p_1=\fc{b(3b-1)}{2}\right),$\
 $\tilde
P_{\nu_{0}}=Z_{p_{0}p_{0}'}$   (by Proposition
5.2.2 of  \cite{MoRoTe}), and $\tilde P_1=\tilde P_{\nu_{0}-b}=1.$

One can compute the order of each $\tilde P_i$ (in terms of the number of
positive
half-twists that appear in the presentation) which is 
$132\; (= 12 \cdot 11)$ 
for a 6-point, and
$12\; (= 3 \cdot 4)$ 
for a 3-point.  
(Since this proof is for the specialist, I shall not give
details of the
calculations).
We sum up the degree of all factors in $\prod\limits_{i=1}^{\nu_0}\tilde
C_i\tilde
P_i,$ and get $p(p-1)$. A priori, $\prod\limits_{i=1}^{\nu_0}\tilde
C_i\tilde P_i$
is part of a braid monodromy factorization of $\Dl_p^2,$ associated to
$S_{1(a,b)}.$
Since the degree of $\Dl_p^2$ is exactly $p(p-1)$ (Proposition 2.3), we get
$\Dl_p^2=\prod\limits_{i=1}^{\nu_0} \tilde C_i\tilde P_i,$ and thus there
are no
extra factors in the braid monodromy factorization; and $\tilde P_i$ is the
local
braid monodromy around
$b_i$.
\hfill
$\qed$\edm
\remark{Remark} The sources of $\tilde C_i$ are the intersection of the lines
$\ell_i$ and $\ell_j,$ for $i$ and $j$, such that $L_i$ and $L_j$ do not
intersect.\endremark

\flushpar{\bf About the computation of $\tilde P_i$ from Theorem 3.2}

Each singular point of $S_0$  is either a 3-point (lies on 3 planes and 2
lines) or a
6-point (lies on 6 planes and 6 lines).  Two intersection lines meet in
each 3-point,
and 6 intersection lines meet in the 6-point. Different types of 3-points,
6-points
are presented in Fig. 3.4.

\midinsert
\centerline{
\epsfysize=2in
\epsfbox{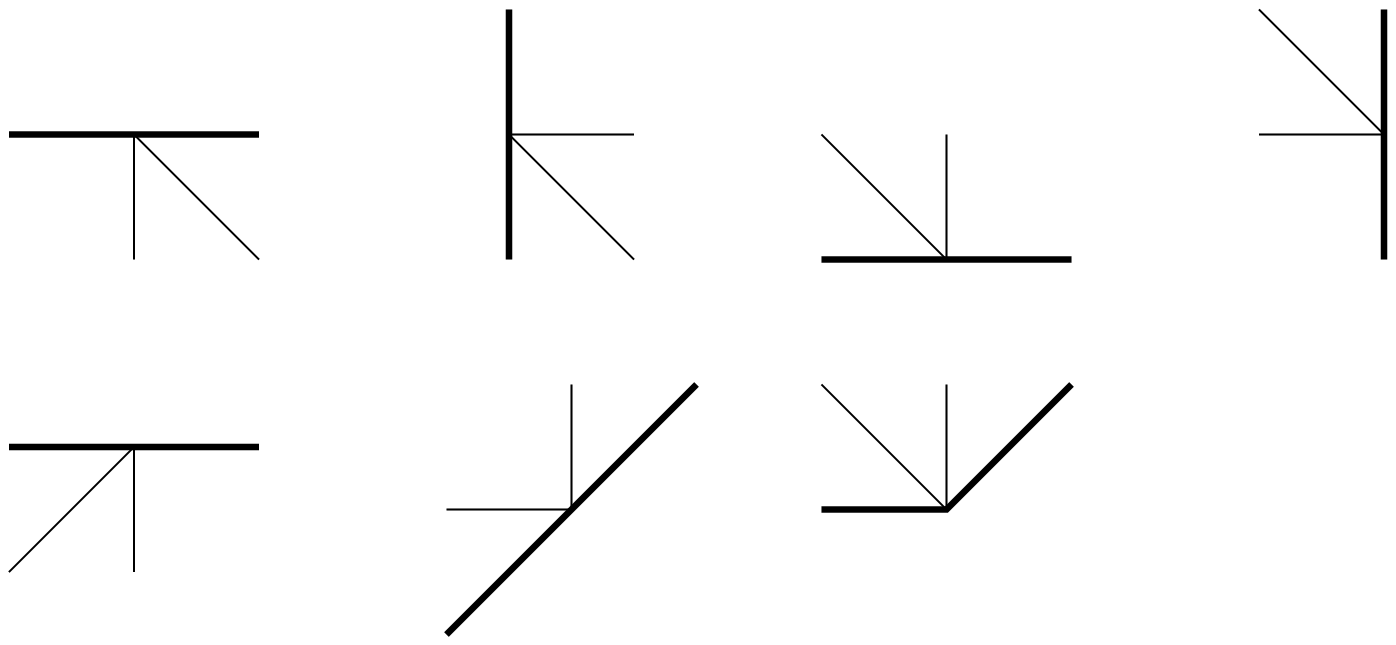}}  
\botcaption{\centerline{3-points}\newline
\centerline{Fig. 3.4(a)}}\endcaption
 \endinsert



\midinsert
\centerline{
\epsfysize=1.2in
\epsfbox{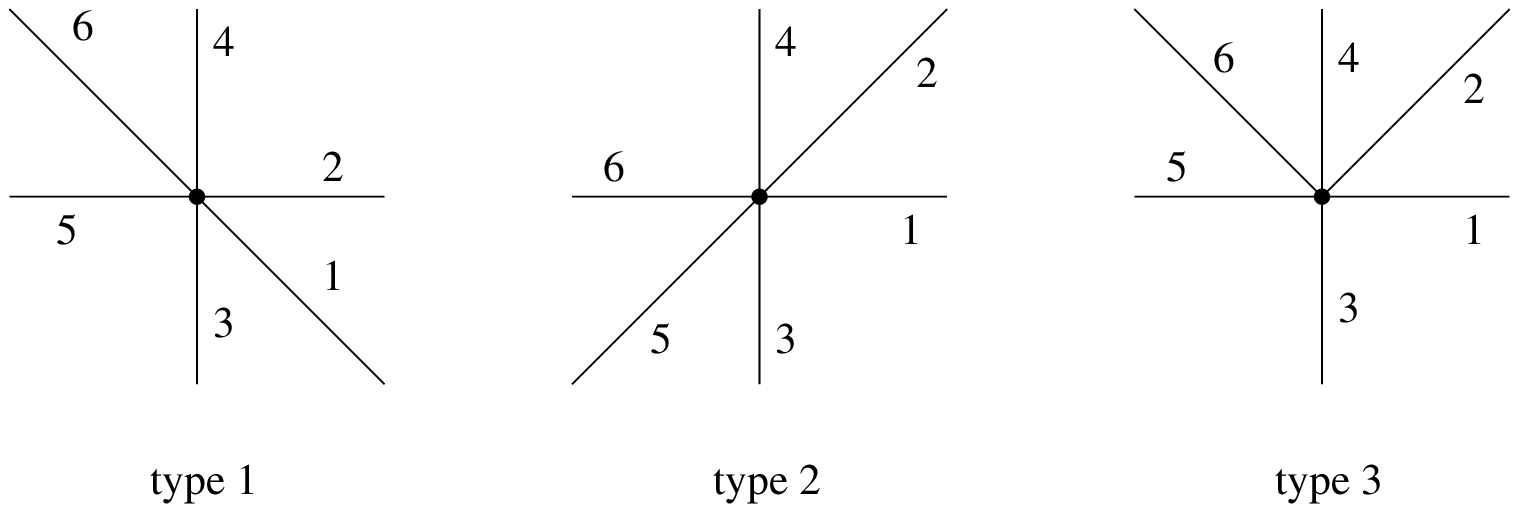}}  
\botcaption{\centerline{6-points}\newline
\centerline{Fig. 3.4(b)}}\endcaption
 \endinsert

(The intersection lines are thinner; the thick line is the border of the
configuration.) ({\it Warning}: In \cite{MoTe1} and \cite{MoTe3}, we refer to
3-points as 2-points, i.e., by the number of lines and not by the number of
planes.)

The difference between the various 3-points lies in the order in which the
lines
appear in the degeneration process.
More precisely, whether the smaller indexed line is a diagonal, vertical, or
a horizontal line and whether the 2 lines meet in the  endpoint with higher 
index of both, or in the endpoint with smaller index of both. 
This difference affects the local braid monodromy
around
each point.

By \cite{MoRoTe}, Prop. 4.4.1 for $a_i$ a 3-point,\ $a_i=L_j\cap L_k$ \ $(L_k$
diagonal) we have
$\tP_i=Z_{k,jj'}^{(3)}\cdot\tilde Z_{kk'}$ where:
$$\tilde Z_{kk'}=\cases  (Z_{kk'})_{Z_{k'j}Z_{jj'}\1Z_{j'k}}&\quad k<j
\\  (Z_{kk'})_{Z_{kj}Z_{jj'}\1Z_{jk'}}&\quad k>j\endcases$$
where $Z_{ij}$ is the half-twist corresponding to a path which connects
$q_i$ with
$q_j$ from below the real line ($q_i$ are real),
and
$(A)_B=B\1AB$ is a conjugation symbol. (Note that a half-twist conjugated by a
half-twist
gives a third half-twist).
$$Z_{k,jj'}^{(3)}=Z_{kj}\cdot Z_{kj'}\cdot(Z_{kj' })_{Z_{jj'}}.$$
is a short notation for product of three half-twists.

Concerning 6-points, we have three types; each one
has six lines meeting in one point.
 We introduce a
local numeration on each configuration which is compatible with the global
ordering.

In \cite{MoTe1,MoTe6}   a complete
computation  for
$\tilde P_i$ where
$a_i$ is a 6-point of type 1 is given
(table
 $\Dl_{\a}^2,$\ $\a$  6-point in \cite{MoTe1} and Lemma 1.1
in
\cite{MoTe6}). In \cite{MoTe1} the local numeration is described by Fig. 3.5:

\midinsert
\centerline{
\epsfysize=1.2in
\epsfbox{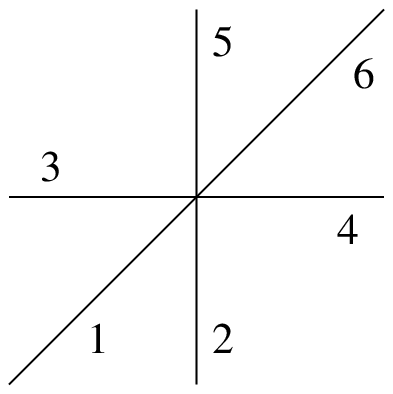}}  
\botcaption{Fig. 3.5}\endcaption
 \endinsert

\flushpar which is obtained by  a $90^\circ$  clockwise turn of the
diagram describing a 6-point of type 1. Thus $\tP_i$ for this type
is determined by the computations given there.

A 6-point of type 2 or 3 is different; no exchange of numeration will
result in the
``classical" 6-point from \cite{MoTe1} with the same order of regeneration.
For these cases, the computations will appear in \cite{AmTe} and \cite{CiMiTe}.

\bk

\head{\S4.\ Fundamental Groups of Galois Covers of Hirzebruch Surfaces}\endhead

 After computing the braid monodromy of a   curve $S$ in $\CP^2$, we can
use the
Zariski-Van Kampen theorem (\cite{VK}, (cf., for example,
\cite{Te1})) to
get a finite presentation for the fundamental group of its complement.
If $S$ is a branch curve of a generic projection from a surface,  the
fundamental group $\pi_1(\CPt\setminus S)$ can yield  a new
invariant of a surface (see \cite{Te3}), using the fact that  in many cases
such groups  are almost polycyclic. Moreover, the fundamental group of the
Galois
cover of the surface is isomorphic to a quotient of a subgroup of
$\pi_1(\CPt\setminus S)$, so we can recover such fundamental groups of
surfaces and
in particular those of the Galois cover of Hirzebruch surfaces.

Let us recall the definition of a Galois cover:

\definition{Definition} \ \underbar{Galois cover w.r.t.
 to generic projection}

Let $X$ be a surface and let $f:X\to\CPt$ be a generic projection of $\deg n.$
Let $ \underset{\underbrace{\ \quad f\qquad
f  \quad}_{n}\quad }\to{X\times\dots\times X} $
 be the fibered product, $$ X\underset f\to\times
\dots
\underset f\to\times X=\{(x_1,\dots, x_n)\bigm|\ \forall i,j,\quad
f(x_i)=f(x_j)\};$$
and let $\Dl$ be the ``big" diagonal, $$\Dl=\{(x_1,\dots,x_n)\bigm|
\exists\ i,j\quad
x_i=x_j\}.$$
Then we define the Galois cover $X_{\Gal}$ of $X$ to be the surface
$$X_{\Gal}=\overline{X\underset f\to\times\dots\underset f\to\times X\setminus
\Dl}.$$ There is a natural projection $\tilde f: X_{\Gal}\to\CPt$ (projection
on the first coordinate).
\enddefinition
The following theorem is concerned with the Galois cover of Hirzebruch surfaces
$F_{k(a,b)}$ (see \S1).
\proclaim{Theorem 4.1} $\pi_1(F_{k(a,b)})_{\Gal}=(\Z_c)^{n-2}$, where
$c=\gcd (a,b),$\ $n=\deg F_{k(a,b)}=2ab+kb^2$.\ep

\demo{Proof} See \cite{MoRoTe} and \cite{FRoTe}.
Here we shall only recall the connection of $\pi_1(X_{\Gal})$ with
$\pi_1(\CPt\setminus S),$ for    $S$ the branch curve of $X\overset f
\to\ri\CPt$
generic.
Let us generically choose  an affine piece $\BC^2$ of
$\CPt.$ Let $X_{\Gal}^{\Aff}$ be the part of $X_{\Gal}$ lying over it.

There is a natural epimorphism $\pi_1(\C^2\setminus S,u_0)\overset
\psi\to\rightarrow S_n$ for   $u_0$ any point not in $S$ and $S_n$ the
symmetric
group on $n=\deg f$
objects. In fact, lifting a loop at $u_0$ to $n$ paths in $X,$ induces a 
permutation of $f\1(u_0).$
Since $\# f\1(u_0)=n,$ we thus get an element of $S_n.$
Clearly, $\psi$ is  surjective.
So we have an 
exact sequence $1\to\ker \psi\to\pi_1(\C^2\setminus S,u_0)\to S_n\to\1$ 
of groups.

In order to establish an isomorphism  of $\pi_1(X_{\Gal}^{\Aff})$ with a
quotient of
a subgroup of $\pi_1(\CPt\setminus S)$,  we have to choose a certain system of
generators for
$\pi_1(\C^2\setminus S,u_0).$
Let $p,u,K,\BC_u^1$ be as in \S2.

 Let $\{\G_j\}_{j=1}^p$ be a $g$-base of $\pi_1(\BC_u^1\setminus K).$
(Recall from \S2
that
$\{\G_j\}$ is a free base).
There is a natural surjection $\pi_1(
{\C}_u^1\setminus S,u)\overset\mu\to\rightarrow
\pi_1(\C^2\setminus S,u)$ induced from the inclusion $\BC_u^1\setminus
S\hookrightarrow \BC^2\setminus S.$
 By abuse of notation, we shall denote the image of $\G_j$  in
$\pi_1(\C^2\setminus S)$ also by $\G_j.$
Clearly, the set $\{\G_j\}_{j=1}^p$ then generates $\pi_1(\C^2\setminus
S,u).$

Since $f$ is stable, the ramification is of order 2 and  $\psi(\G_j)$ is a
transposition in
$S_n.$ So $\G_j^2\in \ker\psi.$
Let  ${\boldsymbol\Gamma}$ be the normal subgroup generated by
$\{\G_j^2\}_{j=1}^p.$
Then ${\boldsymbol\Gamma}\subseteq\ker \psi.$
By the standard isomorphism theorems, we have:
$$1\to \fc{\ker\psi}{{\boldsymbol\Gamma}}\to\fc{\pi_1(\C^2\setminus
S,u)}{{\boldsymbol\Gamma}}\to S_n\to1.$$
In \cite{MoTe1},0.3, we proved $\pi_1(X_{\Gal}^{\Aff})\simeq \fc{\ker
\psi}{{\boldsymbol\Gamma}}.$
In \cite{MoTe8}, we considered the projective case and proved that
$$\pi_1(X_{\Gal})\simeq
\fc{\ker\psi}{\big\la{\boldsymbol\Gamma} ,\prodl_{j=1}^q\G_j\big\ra}$$
This established the connection between $\pi_1(\C^2\setminus S)$ and
$\pi_1(X_{\Gal}).$
The actual 
deduction of $\pi_1(X_{\Gal})$ from  $\pi_1(\C^2\setminus S)$ involves the
Reidemeister-Schreier method from \cite{KMS}.

\hfill
$\qed$\edm
\proclaim{Corollary 4.2} $(F_{k(a,b)})_{\Gal}$ is simply connected iff
$a,b$ are
relatively prime.\ep
\bk

\head{\S5.\ Chern Numbers of Galois Covers of Hirzebruch Surfaces}\endhead

For any    generic (stable, finite) morphism $g: X\to\CPt, $ from a nonsingular
algebraic surface, it can
be shown that the induced
$\XGal$ is nonsingular (see \cite{Te5}). Moreover, if $S\subset \CPt$ is
the branch
curve of
$g,$ and
$\tilde S\subset \XGal$ is the ramification curve of $\tilde g:
X_{\Gal}\to\CPt$ and
$\ell$ a line on
$\CPt$, then the canonical class $K_{\XGal}$ of $\XGal$ is equal to
$g^*(-3\ell)+\tilde S$.
On the level of divisor classes, we have\linebreak $\tilde
S=\fc{1}{2}g^*(S)=\fc{1}{2} mg^*(\ell),$ with  $m=\deg S.$ Hence for the
canonical
class, we get
$K_{\XGal}=\left(\fc{m}{2}-3\right)g^*\ell.$
Thus when $m>6,$ the bundle $K_{\XGal}$ is ample and $\XGal$ is a minimal
surface of
general type.
Moreover, $\XGal$ is a spin manifold iff $K_{\XGal}$ is even iff $m$ is not a
multiple of 4.

\demo{Notation} Let us denote for short $Y_{k(a,b)}=( F_{k(a,b)})_{\Gal}.$\edm

By the above and Lemma 3.1 we get
\proclaim{Corollary 5.1} $Y_{k(a,b)}$ is of general type if and only if one
of the
following is true:
\roster\item"(i)"$ k=0;$\quad
 $ ab\ge 3;$
\item"(ii)"
 $ k=1,2;\quad ab\ge 2;$\item"(iii)"$ k\ge 3.$\endroster

$Y_{k(a,b)}$ is a spin manifold when any of the following is true:
\roster\item"(i)"$ b\equiv 0(4),\qquad  a\equiv 1(2);$
\item"(ii)" $b\equiv 1(4),\qquad k\equiv
0(2);$
\item"(iii)" $b\equiv 2(4), \qquad a+k\equiv 1(2) ;$
\item"(iv)" $ b\equiv 3(4).$\endroster
\ep

\remark{Remarks} (1)\ There are other $Y_{k(a,b)}$ which admit a  spin
structure.

(2)\ Corollary 5.1 is in fact a restatement of Theorem 0.3 b) and c)
from\linebreak
\cite{MoRoTe}.
\endremark

We shall give here a formula from \cite{Te2}
for the Chern numbers of $X_{\Gal}$ in terms of
certain invariants of $X$.

\proclaim{Theorem 5.2} Let $E$ be the  hyperplane section and $K$  the
canonical
divisor of
$X,$ and let $n=\deg X.$ Then
$$  c_1^2(X_{\Gal})=\frac{n!}{4}[(E\cdot K)^2+6n(E\cdot K)+9n^2-12(E\cdot
K)-36n+36],
$$
$$\align c_2&(X_{\Gal})=\\ &\frac{n!}{24}[72-10c_1^2(X)-54(E\cdot
K)-114n+27n^2+14c_2(X)+3(E\cdot K)^2+18n(E\cdot K)].\endalign$$\ep

\demo{Proof}\cite{Te2}, Proposition 2.1.\edm

\proclaim{Proposition 5.3}
For $a,b\ge1$ and  $n=2ab+kb^2(=\deg F_{k(a,b)}),$ we have:
$$\allowdisplaybreaks\align
&c_1^2(Y_{k(a,b)})=\frac{n!}{4}\left\{
\aligned & 4a^2+4b^2-64ab+24a+24b-24a^2b-24ab^2\\&+36+36a^2b^2
+k(12b+4ab-12b^3+36ab^3-24ab^2\\
&-32b^2)+k^2(b^2-6b^3+9b^4)\endaligned\right\}\\
&\ \qquad\qquad =\frac{n!}{4}\left\{
  k^2 b^2 (3b-1)^2 + 4kb(3b-1)(3ab-a-b-3) + 4(3ab-a-b-3)^2\right\},\\
\quad\\
 &c_2(Y_{k(a,b)})=\fc{n!}{8}\left\{
\aligned &
4(4+9a+9b-17ab+a^2+b^2+9a^2b^2-6a^2b-6ab^2)\\
&+2k(9b-17b^2+18ab^3 +2ab-12ab^2-6b^3)\\
&+k^2(9b^4+b^2-6b^3)\endaligned\right\}\\
&\ \qquad\qquad =\frac{n!}{8}\left\{
(3b-1)^2(2a+kb)^2+(9-17b-6b^2)(4a+2kb)+4(b^2+9b+4)\right\}.
\endalign$$
\ep


\demo{Proof} We have here:
$$\align &E(F_{k(a,b)})=aC+bE_0,\\
&K(F_{k(a,b)})=-2E_0+(k-2)C,\\
&c_1^2(F_{k(a,b)})=8,\\
&c_2(F_{k(a,b)})=4,\\
&C\cdot E_0=1,\quad E_0^2=k,\quad C^2=0.\endalign$$
Thus, $$E\cdot K=-2a-2b-bk,$$
$$n=\deg F_{k(a,b)}=E^2=2ab+b^2k.$$
We substitute this in the  formulas from Theorem 5.2 to get the
proposition.\hfill
$\qed$\edm

\proclaim{Corollary 5.4}
For $k=1,$\ $a\ge1,$ $b\ge1$ we have:
$$\align&c_1^2(Y_{1(a,b)})=\fc{(2ab+b^2)!}{4}
(3b^2+6ab-3b-2a-6)^2,\\
&c_2(Y_{1(a,b)})=\fc{(2ab+b^2)!}{8}\left\{\aligned&
16+54b+36a-64ab+4a^2-29b^2+36a^2b^2\\&
-24a^2b-48ab^2-18b^3+9b^4+36ab^3\endaligned\right\}.\endalign$$

For $k=0,$\ $a\ge1,$\ $b\ge1$ we have:
$$\align&c_1^2(Y_{0(a,b)})=(2ab)!(3ab-a-b-3)^2,\\
&c_2(Y_{0(a,b)})=\fc{(2ab)!}{2}  \{
4+9a+9b-17ab+a^2+b^2+9a^2b^2-6a^2b-6ab^2  \}.\endalign$$\ep

Using Theorem 4.1 and Corrolary 5.4, one can get examples of surfaces
with the same
Chern numbers and different fundamental groups:

\proclaim{Theorem 5.5} Let $s,t$ be odd relatively prime positive numbers, then
$$\allowdisplaybreaks\align& \pi_1(Y_{1(s,2t)})=0,\qquad
\pi_1(Y_{0(s+t,2t)})=(\Z_2)^{4st+4t^2-2},\\
\quad\\
&c_1^2(Y_{1(s,2t)})=c_1^2(Y_{0(s+t,2t)})=(4st+4t^2)!
 \left\{\aligned& 9+6s+18t+s^2-30st\\
&-27t^2-12ts^2-48t^2s-30t^3\\
&+36t^2s^2+72st^3+36t^4\endaligned\right\},\\
\quad\\
&c_2(Y_{1(s,2t)})=c_2(Y_{0(s+t,2t)})=\fc{(4st+4t^2)!}{2}
\left\{\aligned& 4+27t+9s-32st\\
&+s^2-29t^2+36s^2t^2-12s^2t\\
&-48st^2-36t^3+36t^4+72st^3\endaligned\right\}.\endalign$$\ep

Using the Hirzebruch formulae for the signature of a
surface in terms of the Chern numbers,
 $\tau(Y)=\frac{1}{3}(c_1^2(Y)-2c_2(Y)),$ our 
expression for the Chern numbers of $Y_{k(a,b)}$ obtained in Proposition 5.3
yields the following result:

\proclaim{Proposition 5.6}
$\tau(Y_{k(a,b)})=\frac{(2ab+kb^2)!}{12}\{4(ab-3a-3b+5)+2k(b-3)b\}.$\ep

In view of the Watershed Conjecture and its role in the "geography of
surfaces", as mentioned in the Introduction, 
the following corollary (0.3 in \cite{MoRoTe}) is interesting:
\proclaim{Corollary 5.7} Let $a\ge 1.$ Then $\tau(Y_{k(a,b)})>0$
if and only if one of the
following
is true:
$$\alignat3
&k=0,\qquad && a\ge 8,\qquad && b=4;\\
& k=0,\qquad && a\ge 6,\qquad &&b\ge 5;\\
&k=1,\qquad && a\ge 6,\qquad && b=4;\\ &k=1,\qquad && a\ge 3, \qquad && b=5;\\
&k=1, \qquad && a\ge 2, \qquad && b=6; \\ & k=1, \qquad && a\ge 1, \qquad
&& b\ge
7;\\ &k=2,\qquad && a\ge 4, \qquad && b=4; \\ & k=2, \qquad && a\ge 1,
\qquad && b\ge
5;\\ &k=3, \qquad && a\ge 2, \qquad && b=4; \\ & k=3, \qquad &&  a\ge 1,
\qquad &&
b\ge 5;\\ &k\ge 4 ,\qquad && a=1, \qquad && b\ge 4.\endalignat$$\ep

{}From 5.1, 5.6 and 5.7 above, we get the following theorem which appeared in a
different form in \cite{MoRoTe}, Theorem 0.4:

\proclaim{Theorem 5.8} Let $a\ge 1.$

\rom{a.} $Y_{k(a,b)}$
is simply connected, of general type, of zero
signature, and a spin manifold if one of the following is true:
\roster\item"(i)"
$k=0, \quad a=7, \quad b=4;$
\item"(ii)" $k=1,\quad a=5, \quad b=4;$
\item"(iii)" $k=2, \quad a=3, \quad b=4;$
\item"(iv)" $k=3, \quad  a=1, \quad b=4.$\endroster

\rom{b.} There are infinitely many triples $k,a,b$, for which $Y_{k(a,b)}$
is
simply connected, of general type, of positive signature, and a spin manifold.
These include infinitely many triples for which $k=1$ or $k\ge 1.$
We present, for example,
$$k=1,\qquad a=3,\qquad b=5,6,7,8, \ldots$$\ep

\bk

\head{\S6.\ Intermediate Galois Covers}\endhead
As mentioned in the introduction, it was wrongly 
conjectured that simply connected surfaces of general type only exist 
in the range $c_{1}^2/c_{2} < 2$ for the ``Chern quotient''. The first 
counterexamples, constructed in 1985, had the Chern quotient 
$c_{1}^2/c_{2}$ just above 2. Those surfaces
were $Y_{0(a,b)}$ for certain choices of $a,b$. 
This was still rather far away from the maximum value 
$c_{1}^2/c_{2} = 3$ that follows from the famous inequality of 
Miyaoka and Yau. Since a surface with $c_{1}^2/c_{2} = 3$ is a free 
quotient of the unit ball in~$\CC^2$, it can never have finite 
fundamental group; in particular, it can never be simply connected. 
It is thus of interest to find out how close one can get to the 
quotient $c_{1}^2/c_{2} = 3$.

 In order to obtain {\it spin} simply connected algebraic surfaces
with positive signature (as in 5.8) with
$ c_1^2/c_2$ closer to 3, we take an intermediate step in the
fibered product. We defined {\it intermediate Galois covers or} $\ell$-{\it th
Galois
cover} as the surface obtained from a fibered product taken $\ell $ times for
$\ell<\deg
X.$ In fact, these constructions give us $ c_1^2/c_2$ closer to 3.

In \cite{Te2}, Theorem 1, we computed the Chern numbers of the $\ell$-th
Galois cover
in terms of $\ell,$\  $\deg X,$ and the following invariants connected to
$S,$ the
branch curve of
$X\to
\CPt:\text{degree}\ (=m),$ number of cusps  $(=\rho),$ number of nodes
$(=d),$ and  $
\deg S^*=\deg  S^{\text{dual}} \ (=\mu).$ For the branch curve of the
Hirzebruch
surface, these invariants were computed in \cite{MoRoTe} and in
\cite{MoTe2}, and the
results are as follows:

\proclaim{Lemma 6.1}
Let $F_{k(a,b)}$ be as above \rom{(\S1)}.

\flushpar\rom{(i)}\ If $a\ge 1,$ then

$n = 2ab+kb^2,$

$m=6ab-2a-2b+k(3b^2-b),$

$\mu=6ab-4a-4b+4+k(3b^2-2b),$

$\vp=24ab-18a-18b+12+k(12b^2-9b),$

\flushpar\rom{(ii)}\ If $a=0,\ k=1,$ then

 $n=b^2,$

$m=3b(b-1),$

$\mu=3(b-1)^2,$

$\vp=3(b-1)(4b-5),$

$d=\df{3}{2}\ (b-1) (3b^3-3b^2-14b+16).$
\ep

\demo{Proof}

(i) \cite{MoRoTe}, Lemma 7.1.3.

(ii) \cite{MoTe2}, \S2.\hfill$\qed$\edm

If one substitutes the above $n,\ m,\ \mu,\ d,\ \vp$ in the formulas of
\cite{Te2},
  one gets  the Chern classes of the $\ell$-th Galois
cover of
$F_{k(a,b)}.$

Theorem 6.2 treats the case $\ell=n$ (the full Galois cover), and Corollary 6.3
treats an intermediate step $\ell=4<9=n.$

\proclaim{Theorem 6.2}
For $k=1,$\ $a=0$ \rom{(}Veronese surface of order $b$\rom{)} we get
$$c_1^2(Y_{1(0,b)})=\fc{9}{4}(b^2)!\{ b^4-2b^3-3b^2+4b+4\},$$
$$c_2(Y_{1(0,b)})=\fc{(b^2)!}{8}\{16+54b-29b^2-18b^3+9b^4\}.$$\ep

\demo{Proof} \cite{Te2}, Theorem 1 and Lemma 6.1(b).\hfill$\qed$\edm

Computing $c_1^2/c_2$ for $\ell=n$  gives us numbers close to the line
$c_1^2=2c_2.$
On the other hand:
\proclaim{Corollary 6.3} For $k=1,\ a=0,\ b=3$ and $\ell=4,$  we get
$ c_1^2/c_2=2.73$ (almost precisely).\ep
\remark{Remark}
By experimental substitution it seems that for large $b,$ the signature
$=\fc{1}{3}(c_1^2-2c_2)$
changes from negative to positive around  $\ell=\fc{3}{4}n.$ \endremark

\enddocument